\documentclass{article}
\usepackage{amsmath,array,theorem,doublespace}
\usepackage{tabularx,xspace}
\usepackage{commands,shell}
\ifpdf
\usepackage[pdftex]{graphicx}
\else
\usepackage[dvips]{epsfig}
\fi
\usepackage{amssymb}
\begin{document}

\title{Output Feedback Control for Stabilizable and Incompletely
Observable Nonlinear Systems\thanks{This work 
was supported by NASA Glenn Research Center, Grant NAG3-2084.}}
\author{\bf
Manfredi Maggiore\hspace{1cm} 
Kevin Passino \\
\\
\it
Department of Electrical Engineering,
The Ohio State University \\
\it 
 2015 Neil Avenue, Columbus, OH
43210-1272 \\
}
\date{}

\maketitle
\begin{abstract}
This paper introduces a new approach for  output
feedback stabilization of SISO systems which, unlike most of the techniques
found in the literature, does not use high-gain observers and control
input saturation to achieve separation between the state feedback  and
observer designs. Rather, we show that by using nonlinear observers,
together with a projection algorithm, the same kind of separation
principle is achieved for a larger class of systems, namely
stabilizable and incompletely observable plants. Furthermore, this new
approach avoids using knowledge of the inverse of the observability
mapping, which is needed by most techniques in the literature  when controlling general
stabilizable systems.

\end{abstract}
\begin{spacing}{1.3}
\vspace*{-1cm}
\section{Introduction}
 The area of  nonlinear output feedback control has
received much attention after the publication of the work 
\cite{EsfKha92}, in which the authors developed a systematic
strategy for the output feedback control of input-output linearizable systems
with full relative degree, which employed two basic tools: an high-gain
observer to estimate the derivatives of the outputs (and hence the
system states in transformed coordinates), and control input
saturation to isolate the peaking phenomenon of the observer from the
system states. Essentially the same  approach has  later been applied in  a number of papers 
by various  researchers (see, e.g., \cite{KhaEsf93,LinSab95,Jan96,MahKha96,
MahKha97,Isi97,AtaKha99}) to solve different problems in
output feedback control. In most of the papers found in the literature,
(see, e.g., \cite{EsfKha92,KhaEsf93,LinSab95,MahKha96,Jan96,MahKha97})
the authors consider input-output \FLE systems with either full
relative degree or minimum phase zero dynamics. 
The work in \cite{Tor92} showed that for nonminimum phase systems the
problem can be solved by extending the system with a  chain of
integrators at the input side. However,  the results contained there
are local.
In \cite{TeePra94}, by putting together this idea 
with the approach found in \cite{EsfKha92},
  the authors were able to show how to
solve the  output
feedback stabilization problem for general stabilizable and uniformly
completely  observable systems, provided that the inverse of
the observability mapping is explicitly known. The recent  work in \cite{AtaKha99}
unifies all these approaches  to prove  a separation principle 
 for a very general class of nonlinear systems.
It appears that the largest class of nonlinear SISO systems for which 
the output feedback stabilization problem has been solved is that
of locally stabilizable and completely observable systems. 
Moreover, when dealing with systems which are not feedback
linearizable, the works \cite{Tor92,TeePra94, AtaKha99} require the explicit
knowledge of the inverse of the observability mapping, thus somewhat
restricting the variety of problems to which their algorithm can be
applied.

The objective of this paper is to relax the two restrictions above,
by developing a new  output feedback strategy for nonlinear SISO
locally or globally stabilizable systems which are only observable on 
regions of the state
space. Furthermore,  for the implementation of our controller,
the inverse of the observability mapping is not needed.
These two features are achieved by means  of a nonlinear
observer instead of the standard high-gain observer found in the
literature, and of a new projection algorithm which {\em eliminates}
the peaking phenomenon in the observer states, thus avoiding the
need to use control input saturation. To the best of our knowledge,
this work, besides introducing  a new methodology for 
output feedback control design, enlarges the class of SISO systems
considered in the literature of the field so far. 

\section{Problem Formulation and Assumptions}
Consider the following dynamical system,
\begin{equation}
\begin{split}
&\dot x = f (x,u)  \\
&y = h(x,u) \label{eq:system}
\end{split}
\end{equation}
where $x \in \Re^n, u, y\in \Re$, $f$ and $h$ are known smooth
functions, and $f(0,0) = 0$.
Our  control objective is to construct a stabilizing controller 
for (\ref{eq:system}) without the availability of the system states
$x$. In order to do so, we need an observability assumption. Define 
\begin{equation}
y_e \eqdef  \begin{bmatrix} 
y \\
\vdots \\
y^{(n-1)}
\end{bmatrix} = 
\HH\left(x, u, \ldots, u^{(n_u-1)}\right) \eqdef \left[
\begin{array}{c}
h(x,u) \\
\varphi_1(x,u,u^{(1)}) \\
\vdots \\
\varphi_{n-1} \left(x,u, \ldots, u^{(n_u-1)}\right)
\end{array}\right]
\end{equation}
($y^{(n-1)}$ is the $n-1$-th derivative) where 
\begin{equation}
\begin{split}
&\varphi_1(x,u,u^{(1)})= \DER{h}{x} f(x,u) + \DER{h}{u} u^{(1)}\\
&\varphi_2(x,u,u^{(1)},u^{(2)}) = \DER{\varphi_1}{x} f(x,u) +
\DER{\varphi_1}{u} u^{(1)} + \DER{\varphi_1}{u^{(1)}} u^{(2)}\\
& \hspace{2cm}\vdots \\
&\varphi_{n-1} \left(x,u, \ldots, u^{(n_u-1)}\right) = \DER{\varphi_{n
- 2 }}{x} f(x,u) + \DER{\varphi_{n-2}}{u} u^{(1)} + \ldots +
\DER{\varphi_{n-2}}{u^{(n_u-2)}} u^{(n_u-1)}
\end{split}
\label{eq:phi_defn}
\end{equation}
where $0 \leq n_u \leq n$ ($n_u=0$ indicates that there is no
dependence on $u$). In the most general case, $\varphi_i =
\varphi_i(x, u, \ldots, u^{(i)})$, $i=1, 2, \ldots, n-1$. In some cases, however, we may
have that $\varphi_i = \varphi_i(x,u)$ for all $i = 1 , \ldots , r-1$
and some integer $r > 1$. This happens in particular when system 
(\ref{eq:system}) has a well-defined relative degree $r$. Here, we
do not require the system to be input-output \FLE, and hence to
possess a well-defined relative degree. In the
case of systems with well-defined  relative degree, $n_u=0$
corresponds to having $r \geq n$, while $n_u = n$
corresponds to  having   $r=0$.
Now, we are ready to state our first assumption.
\ass{1}. System (\ref{eq:system}) is observable over the
set $\X \times \U \subset \Re^n \times \Re^{n_u}$ containing the origin,
\ie  the mapping 
\begin{equation}
y_e = \HH\left(x,u, \ldots, u^{(n_u-1)}\right) 
\label{eq:observability}
\end{equation}
is invertible with respect to $x$ and its inverse is smooth, for all
$x \in \X$, $[ u , u^{(1)}, \ldots , u^{(n_u-1)}]^\top \in \U$.
\begin{rem} In the existing literature, an assumption similar  to A1 can
be found in \cite{Tor92} and \cite{TeePra94}. 
It is worth stressing, however, that in that work  the authors adopt
a global observability assumption, \ie the set $\X \times \U$ is taken to be
$\Re^n \times \Re^{n_u}$. In many practical applications the system under consideration
may be observable in some subset of $\Re^n \times \Re^{n_u}$ only, thus preventing the
use of most of the output feedback techniques found in the
literature, including the ones found in  \cite{Tor92}, \cite{TeePra94}, and
\cite{AtaKha99}. 
\end{rem}
Next, augment the system with $n_u$ integrators on the input side,
which corresponds to using a compensator of order $n_u$.
System (\ref{eq:system}) can be rewritten as follows,
\begin{equation}
\begin{split}
&\dot x = f(x,z_1) \\
&\dot z_1 = z_2 \\
&\s \vdots \\
&\dot z_{n_u} = v
\label{eq:sys+integrators}
\end{split}
\end{equation}
Notice that the differential equation (\ref{eq:sys+integrators}) is now 
affine in the new control $v$. Define  the  extended state variable $\chi =
[x^\top , z^\top]^\top \in \Re^{n+n_u} $, and  the associated {\em extended
system} 
\begin{equation}
\begin{split}
&\dot \chi = f_e (\chi) + g_e v \\
& y = h_e (\chi)
\end{split}
\label{eq:extended} 
\end{equation}
where $f_e(\chi) = [ f^\top(x,z_1) ,\,\,\, z_2 \,, \ldots ,\, z_{n_u} ,\, 0 ]^\top$, 
$g_e = [0, \, \ldots, \, 1]^\top$, and $h_e (\chi) = h(x,z_1)$.
\ass{2}. The origin of (\ref{eq:system}) is locally stabilizable
(stabilizable) by a static function of $x$, \ie there
exists a smooth function $\bar u(x)$ such that the origin is 
an asymptotically stable (globally asymptotically stable) equilibrium
point  of $\dot x = f(x,\bar u(x))$.
\begin{rem} Assumption A2 implies that the origin of the extended system
(\ref{eq:extended}) is locally stabilizable (stabilizable) by a
function of $\chi$ as well. A proof of the local stabilizability
property for (\ref{eq:extended}) may be found, e.g., in \cite{Son89},
while its global counterpart is a well known consequence of the
integrator backstepping lemma (see, e.g., Theorem 9.2.3 in
\cite{Isi95} or Corollary 2.10 in \cite{KrsKanKok95}).
Therefore we conclude that for the
extended system (\ref{eq:extended}) there exists a smooth control
$\bar v(\chi)$ such that its origin is asymptotically stable under
closed-loop  control. Let $\DD$ be the domain of attraction of the origin of
(\ref{eq:extended}), and notice that, when A2 holds globally, $\DD =
\Re^n \times \Re^{(n_u)}$.  
\label{rem:stabilizability}
\end{rem}
\begin{rem} 
In \cite{Tor92} the
authors consider affine systems and use a \FLY assumption in place of
our A2. Here,  we consider the more general class of non-affine
systems for which the origin is locally stabilizable (stabilizable). In this respect, our
assumption A2 relaxes also the stabilizability assumption found in
\cite{TeePra94}, while it is essentially equivalent to Assumption 2 
in \cite{AtaKha99}.
In conclusion, the class of systems satisfying A1 and A2 is, to
the best of our knowledge, the largest  class of SISO 
systems considered in the output feedback literature so far.
\end{rem}
\section{Nonlinear Observer: Its Need and Stability Analysis}
Assumption A2  allows
us to design a stabilizing state feedback control
$v = \phi(x,z)$.
In order to perform output feedback control $x$ should be
replaced by its estimate. Many researchers adopted an input-output
\FLY assumption  (\cite{EsfKha92}, \cite{KhaEsf93},
\cite{MahKha96}, \cite{Jan96},  \cite{MahKha97}) and  transformed the system into  normal
form
\begin{equation}\begin{split}
&\dot \pi_i = \pi_{i+1} , \s 1 \leq i \leq r-1 \\
&\dot \pi_r = \bar f (\pi, \Pi) + \bar g (\pi , \Pi ) \, u \\
&\dot \Pi = \Phi (\pi,\Pi) , \s \Pi \in \Re^{n-r} \\
&y = \pi_1
\label{eq:normal_form}
\end{split}\end{equation}
In this framework the problem of output feedback control finds
a very natural formulation, as the first $r$ derivatives of
$y$ are equal to the states of the $\pi$-subsystem (\ie the linear
subsystem). The works \cite{EsfKha92,KhaEsf93,Jan96} solve the output
feedback control problem
for systems with no zero dynamics (\ie $r=n$), so that the first 
$n-1$ derivatives of $y$ provide the entire state of the system.
In the presence of zero dynamics ($\Pi$-subsystem), the use of 
 input-output \FL to put the system into normal form
(\ref{eq:normal_form})
forces the use of a minimum phase assumption (e.g, \cite{MahKha96})
since the states of the $\Pi$-subsystem cannot be estimated from the
derivatives of the output and, hence, cannot be controlled by output
feedback.
It is for this reason that output feedback control of 
nonminimum phase systems is regarded as a particularly challenging
problem. Researchers who have addressed this problem (e.g.,
\cite{Tor92}, \cite{TeePra94}) relied on the explicit knowledge of
the inverse of the mapping $\HH$ in (\ref{eq:observability})
\[
x = \HH^{-1} \left( y_e , z_1 , \ldots , z_{n_u} \right) 
\]
so that estimation of the first $n-1$ derivatives of $y$ (the vector
$y_e$) provides an estimate of $x$
\[
\hat x = \HH^{-1} \left( \hat y_e , z_1 , \ldots , z_{n_u} \right)
\]
since the vector $z$, being the state of the controller, is known.
Next, to estimate the derivatives of $y$, they employed an high-gain
observer.
 Both the works \cite{Tor92} and  \cite{TeePra94} (the
latter dealing with  the larger class of stabilizable systems) rely on
the knowledge of $\HH^{-1}$ to prove closed-loop stability.
In addition to this, the recent work \cite{AtaKha99} proves that a
separation principle holds for a quite general class of nonlinear
systems which  includes (\ref{eq:system}) provided that $\HH^{-1}$ is
explicitly known and that the system is  uniformly completely observable.
In order to develop a practical  output feedback control algorithm,
however, $\HH^{-1}$ cannot be assumed to be 
explicitly known. Hence, rather than designing an high-gain observer
to estimate $y_e$ and using $\HH^{-1} (\cdot, \cdot)$ to get $x$, the
approach adopted here is 
 to estimate $x$ directly using a nonlinear observer for 
system (\ref{eq:system}) and
 using the fact that 
 the $z$-states are known. In
other words, we can regard our problem as that of building a 
reduced order observer for the closed-loop system
states\footnote{Throughout this section we assume A1 to hold globally,
since we are interested in the ideal convergence
properties of the state estimates. In the next section we will show
how to modify the observer equation in order to achieve the same
convergence properties when A1 holds over the set $\X \times \U
\subset \Re^n \times \Re^{n_u}$.}.
The observer  has the form
\begin{equation}
\begin{split}
\dot{\hat x} = f(\hat x,z_1) + &\left[ \DER{\HH(\hat x, z)}{\hat x}
\right]^{-1} \E^{-1}L\, [ y(t) - \hat y(t) ]
 \\
&\hat y(t) = h (\hat x, z_1)
\end{split}
\label{eq:observer}
\end{equation}
where $L$ is a $n \times 1$ vector, $\E = \text{diag} \left[ \rho,
\rho^2, \ldots , \rho^n \right]$, and $\rho \in (0, 1]$ is a
fixed design constant.

Notice that (\ref{eq:observer}) does not require any knowledge of
$\HH^{-1}$ and has the advantage of operating in $x$-coordinates.
The observability assumption A1 implies that the Jacobian of the
mapping $\HH$ with respect to $x$ is invertible, and hence the inverse of
$\partial \HH(\hat x, z) / \partial \hat x$ in (\ref{eq:observer}) is well defined.
In the work \cite{CicDalGer93}, the authors used  an observer
structurally identical to (\ref{eq:observer}), for the more
restrictive class of input-output \FLE systems with full relative
degree. Here, by modifying the definition of the mapping $\HH$, we
considerably relax these conditions by just requiring the general
observability assumption A1 to hold. Furthermore, we propose a
different proof than the one found in \cite{CicDalGer93} which,
besides being easier (in our view), clarifies the relationship among
(\ref{eq:observer}) and the high-gain observers commonly found in the
output feedback literature. Next, we state the result and its proof.
\begin{thm} Consider system (\ref{eq:sys+integrators}) and assume  A1 is
satisfied for $\X = \Re^n$ and  $\U = \Re^{n_u}$, the state $\chi$
belongs to a compact invariant set $\Omega$, and that  $|v(t)|\leq M$ for all
$t \geq 0$, with $M$ a positive constant. Choose $L$ such that $A_c -
L C_c$, where  $(A_c , B_c , C_c)$ is the controllable/observable
canonical realization, is Hurwitz.

%
Under these  conditions and using observer (\ref{eq:observer}), the following two 
properties hold
\begin{itemize}
\item[(i)] Asymptotic stability of the estimation error: 
There exists $\bar \rho$, $0 < \bar \rho \leq 1$, such that for all
$\rho \in ( 0, \bar \rho )$, $\hat x \rightarrow x$ as $t \rightarrow +\infty$.
\item[(ii)] Arbitrarily fast rate of convergence: For each positive $T
, \epsilon$,  there exists $\rho^*$, $0 < \rho^* \leq 1$,
such that for all $\rho \in (0,  \rho^*]$, $\|\hat x - x \|\leq \epsilon \,\, \forall t \geq T$.
\end{itemize}
\label{thm:observer}
\end{thm}
\pf Consider the filtered transformation
\begin{equation}
\xi = \HH(x,z) = \left[
\begin{array}{c}
h(x,z_1) \\
\varphi_1(x,z_1,z_2) \\
\vdots \\
\varphi_{n-1} \left(x,z_1, \ldots, z_{n_u} \right)
\end{array}\right]
\label{eq:tsf:xi}
\end{equation}
A1  guarantees that $x = \HH^{-1}(\xi, z)$ is well-defined,
unique, and smooth. Let us express system (\ref{eq:sys+integrators}) in  new
coordinates. By definition, $\xi = [ y ,  \, \dot y , \, \ldots
, \, y^{(n-1)} ]^\top$ and, with $\varphi_{n-1}$ defined in (\ref{eq:phi_defn}),
\begin{align*}
y^{(n)} &=  \left[ \DER{\varphi_{n-1}}{x} f \left(\HH^{-1}(\xi,z),z\right) 
+ \sum_{k=1}^{n_u-1}
\DER{\varphi_{n-1}}{z_k}\left(\HH^{-1}(\xi,z),z\right) \,
z_{k+1} \right] + \left[
\DER{\varphi_{n-1}}{z_{n_u}}\left(\HH^{-1}(\xi,z),z\right)\right] v \\
&\eqdef \alpha(\xi,z) + \beta(\xi,z)  v
\end{align*}
Hence, in the new coordinates (\ref{eq:sys+integrators}) becomes
\begin{equation}
\dot{\xi} = A_c \xi + B_c \left[
\alpha(\xi,z) + \beta(\xi, z) v \right]
\label{eq:system_transformed}
\end{equation}
Next, transform the observer (\ref{eq:observer}) to
new coordinates  $\hat \xi = \left[ \hat y, \dot{\hat y}, \ldots, \hat
y^{(n-1)} \right]^\top =\HH(\hat x,z)$ so that

\begin{align}
\dot{\hat \xi}_1 &= \DER{h}{\hat x} f(\hat x, z_1) + \DER{h}{\hat x}
\left[ 
\DER{\HH}{ \hat x} \right]^{-1} \E^{-1}L\, [ y - h(\hat x,z_1)]
+ \DER{h}{z_1} \dot z_1 \nonumber \\
&= \hat \xi_2  + \DER{h}{\hat x}  \left[
\DER{\HH}{\hat x} \right]^{-1} \E^{-1}L\, [ y - h(\hat x,z_1)] 
\label{eq:xi_1_hat}
\end{align}
Similarly, for $i = 2, \ldots, n-1$
\begin{equation}
\dot{\hat \xi_i} = \DER{\varphi_{i-1}}{\hat x} (\hat x, z_1, \ldots,
z_{i}) \left\{ f(\hat x,z_1) + \left[ \DER{\HH}{\hat x} \right]^{-1}
\E^{-1}L\, (y-h(\hat x,z_1)) \right\} + \sum_{k=1}^{i}
\DER{\varphi_{i-1}}{z_k} z_{k+1} 
\end{equation}
By definition, 
\[
{\hat \xi}_{i+1} = \varphi_i (\hat x, z_1, \ldots,
z_{i+1}) = \DER{\varphi_{i-1}}{\hat x} f(\hat x, z_1) +
\sum_{k=1}^{i} \DER{\varphi_{i-1}}{z_k} z_{k+1} 
\]
Hence, we conclude that 
\begin{equation}
\dot{\hat \xi_i} = \hat \xi_{i+1} +  \DER{\varphi_{i-1}}{\hat x}
\left[ \DER{\HH}{\hat x} \right]^{-1} \E^{-1}L\, (y-h(\hat x,z_1)),
\s i=2, \ldots , n-1
\label{eq:xi_i_hat}
\end{equation}
Finally, 
\begin{equation}
\dot{\hat \xi}_n = \alpha(\hat \xi,z) + \beta(\hat \xi,z) v +
\DER{\varphi_{n-1}}{\hat x} \left[
\DER{\HH}{\hat x} \right]^{-1} \E^{-1}L\, [ y - h(\hat x,z_1)]
\label{eq:xi_n_hat}
\end{equation}
By using (\ref{eq:xi_1_hat}),
(\ref{eq:xi_i_hat}),
and (\ref{eq:xi_n_hat}) we can write, in compact form,
\begin{align}
\dot{\hat \xi} &= A_c \hat \xi + B_c [\alpha(\hat \xi,z)+
\beta(\hat \xi,z) v ] + \left[
\DER{\HH}{\hat x} \right] \left[
\DER{\HH}{\hat x} \right]^{-1}  \E^{-1}L\, [ y - h(\hat x,z_1)] \nonumber
 \\
&= A_c \hat \xi + B_c [\alpha(\hat \xi,z)+
\beta(\hat \xi,z) v ] + \E^{-1}L\, [ \xi_1 - \hat{\xi}_1]
\label{eq:observer_transformed}
\end{align}
Define the observer error in the new coordinates,
$\tilde \xi = \hat \xi - \xi$. Then, the observer error
dynamics are given by
\begin{equation}
\dot{\tilde \xi} = (A_c - \E^{-1}L\, C_c) \tilde \xi + B_c \left[
\alpha(\hat \xi,z) + \beta(\hat \xi,z) v -
\alpha( \xi,z) - \beta(\xi,z) v \right] 
\label{eq:error_dynamics}
\end{equation}
Next, define the coordinate transformation
\begin{equation}
\tilde \nu = \E^\prime \tilde \xi, \s  \E^\prime
\eqdef \text{diag} \left[ \frac{1}{\rho^{n-1}}, \frac{1}{\rho^{n-2}}, \ldots, 1 \right]
\label{eq:tsf:nu}
\end{equation}
In the new domain the observer error equation becomes
\begin{equation}
\dot{\tilde \nu} = \frac{1}{\rho}(A_c - L C_c) \tilde \nu + B_c
\left[ \alpha (\hat \xi,z) + \beta(\hat \xi,z) v -   \alpha (\xi,z) -
\beta(\xi,z) v \right]
\label{eq:observer:nu}
\end{equation}
where, by assumption,  $A_c - L C_c$ is Hurwitz. Let $P$ be the
solution to the Lyapunov equation
\begin{equation}
P (A_c - L C_c) + (A_c - L C_c)^\top P = -I
\label{eq:lyap}
\end{equation}
and consider the Lyapunov function candidate $V_o(\tilde \nu) = \tilde \nu^\top P
\tilde \nu$. 
Calculate the time derivative of $V_o$ along the $\tilde \nu$
trajectories:
\begin{equation}
\dot V_o = -\frac{\tilde \nu^\top \tilde \nu}{\rho} + 2 \tilde
\nu^\top P B_c \left[ \alpha (\hat \xi,z) + \beta(\hat \xi,z) v -   \alpha (\xi,z) -
\beta(\xi,z) v \right]
\label{eq:W1_dot}
\end{equation}
Next, we seek to provide a bound to the bracketed term in
(\ref{eq:W1_dot}).
Without loss of generality let $[\hat x(0)^\top, z(0)^\top]^\top$ $\in
\Omega$ and define the compact set
$\K_{\tilde \xi} \eqdef \left\{ \tilde \xi \in \Re^n \,|\, \hat
\xi, \xi \in \HH(\Omega) \right\}$. By definition, $ \K_{\tilde \xi}$ contains the
initial condition $\tilde \xi(0)$.
Next, define the set  
$\Lambda_\kappa \eqdef \left\{ \tilde \xi \in \Re^n \,|\, V_o(\E^\prime
\tilde \xi) \leq \kappa \right\}$ where $\kappa$ is chosen so that
$\K_{\tilde \xi} \subset \Lambda_\kappa$. In the following we will
prove that $\Lambda_\kappa$ is
invariant under (\ref{eq:error_dynamics}), and that all trajectories 
originating in $\Lambda_\kappa$ converge asymptotically to the origin.
Recalling that
$\hat \xi = \xi + \tilde \xi$, let $\K_{\hat \xi} \eqdef \left\{
 \hat \xi \in \Re^n \,|\, \xi \in \HH(\Omega) , \tilde \xi \in
\Lambda_\kappa \right\}$, and notice that $\HH(\Omega) \subset \K_{\hat \xi}$.

Due to the smoothness  of $\alpha + \beta v$ and the
boundedness of $v$, the following inequality holds true over the
compact set $\K_{\hat \xi}$:
\begin{equation}
\begin{split}
&\sup_{| v | \leq M}\left\| \alpha(\xi^1,z) + \beta(\xi^1,z) v -
\alpha(\xi^2,z) - \beta(\xi^2,z) v 
\right\| \leq  \gamma \|\xi^1 - \xi^2\| \,, \\  & \forall \xi^1 ,
\xi^2 \in \K_{\hat \xi} , 
\forall z \in \Omega^z \eqdef \{ z \in
\Re^{n_u} \,|\, \chi \in \Omega \}
\label{eq:lipschitz}
\end{split}
\end{equation}
Using this inequality in (\ref{eq:W1_dot})
\begin{equation}
\dot V_o \leq -\frac{\| \tilde \nu \|^2}{\rho} + 2 \|P\| \gamma
\|\tilde \xi\| \|\tilde \nu\| \leq -\frac{\| \tilde \nu \|^2}{\rho} + 2 \|P\|
\gamma \|\tilde \nu\|^2
\label{eq:V_o_dot}
\end{equation}
Defining $\bar \rho = \min\{1/(2 \|P\| \gamma) , 1\}$, we conclude that
for all $\rho < \bar \rho$ the $\tilde \xi$ trajectories starting in $\K_{\tilde
\xi}$ will converge asymptotically to the origin, and hence part (i)
of the theorem is proved.

As for part (ii), note that $\lm(\E' P \E') \geq \lm(\E')^2 \lmP =
\lmP$, since $\lm(\E') =1$. Next, $\lM(\E' P \E') \leq \lM(\E')^2 \lMP =
1 /(\rho^{2(n-1)}) \lMP$, since $\lM(\E') = 1 / \rho^{(n-1)}$. Therefore
\begin{equation}
\lmP \| \tilde \xi \|^2 \leq V_o = \tilde \xi^\top
\E^\prime P \E^\prime \tilde \xi \leq \frac{1}{\rho^{2(n-1)}} \lMP \|
\tilde \xi \|^2
\label{eq:W1_ineq}
\end{equation}
Define $\bar \epsilon$ so that $\| \tilde \xi \| \leq \bar \epsilon$
implies that $\|\hat x -  x\| \leq  \epsilon$ (the smoothness of
$\HH^{-1}$ guarantees that $\bar \epsilon$ is well defined). By
inequality (\ref{eq:W1_ineq}) we have that $V_o \leq \bar \epsilon^2 \lmP$
implies that 
$\| \tilde \xi \| \leq \bar \epsilon$, and $V_o(0) \eqdef V_o
(\tilde \nu (0)) \leq  ( 1 / \rho^{2(n-1)} ) \lMP \|
\tilde \xi(0) \|^2$. Moreover, from (\ref{eq:V_o_dot})
\begin{equation}
\dot V_o(t) \leq - \left( \frac{1}{\rho} - 2 \|P\| \gamma \right) \|
\tilde \nu \|^2 \leq - \frac{1}{\lmP} \left(\frac{1}{\rho} - 2 \|P\|
\gamma \right) V_o(t) 
\end{equation}
Therefore, by the Comparison Lemma (see, e.g., \cite{Yos66}), $V_o(t)$ satisfies
the following inequality
\begin{equation}
\begin{split}
V_o(t) &\leq V_o(0) \exp\left\{- \frac{1}{\lmP} \left(\frac{1}{\rho} - 2 \|P\|
\gamma \right) t \right\} \\
&\leq \frac{1}{\rho^{2(n-1)}} \lMP \|
\tilde \xi(0) \|^2 \exp\left\{- \frac{1}{\lmP} \left(\frac{1}{\rho} - 2 \|P\|
\gamma \right) t \right\} 
\label{eq:W1(t)}
\end{split}
\end{equation}
which, for sufficiently small $\rho$, can be written as 
\[
V_o(t) \leq \frac{a_1}{\rho^{2n}} \exp\left\{- \frac{a_2}{\rho} t
\right\}, \s a_1 , a_2 > 0 
\]
An upper estimate of the time $T$ such that 
$\| \hat x -  x \| \leq \epsilon$ for all $t \geq T$, is calculated as follows 
\[
\frac{a_1}{\rho^{2n}} \exp \left\{ -\frac{a_2}{\rho} t \right\} \leq
\bar \epsilon^2 \lmP \text{ for all } t \geq T = \frac{2n \rho}{a_2}
\log \left( \frac{a_1}{\bar \epsilon \rho} \right)
\]
Noticing that $T \tends  0$ as $\rho \tends 0$, we conclude that $T$
can be made arbitrarily small by choosing a sufficiently small
$\rho^*$, thus concluding the proof of part (ii). 
\BBOX
\begin{rem} Part (ii) of Theorem \ref{thm:observer} implies that the observer convergence
rate can be made arbitrarily fast. This property is essential for 
closed-loop stability.
\end{rem}
\begin{rem} Using inequality (\ref{eq:W1(t)}), we
find the upper bound for the estimation error in $\xi$-coordinates
\begin{equation}
\| \tilde \xi \| \leq \sqrt{\frac{\lMP}{\lmP}} \frac{1}{\rho^{n-1}}
\|\tilde \xi (0) \| \exp \left\{  - \frac{1}{2 \lmP} \left(\frac{1}{\rho} - 2 \|P\|
\gamma \right) t \right\}
\label{eq:observer_rate}
\end{equation}
Hence, during the initial transient, 
$\tilde \xi(t)$ may exhibit peaking, and the size of the peak  grows larger 
 as $\rho$ decreases and the convergence rate is made faster.
 This phenomenon and its
implications on output feedback control has been studied in the
seminal work \cite{EsfKha92}. The analysis  in that paper
shows that a way to isolate the peaking of the observer estimates from
the system states is to saturate the control input outside of the
compact set of interest. 
The same idea has  then been adopted in several other works in the
output feedback control literature (see, e.g.,
\cite{EsfKha92,KhaEsf93,TeePra94,LinSab95,MahKha96,Jan96,MahKha97,Isi97,
AtaKha99}). Rather than following this approach, in the
next section we
will present a new technique for isolating the peaking phenomenon
which allows for the use of the weaker Assumption A1.
\end{rem}
\begin{rem} It is interesting to note that in $\xi$-coordinates the
nonlinear observer (\ref{eq:observer}) is identical to the standard
high-gain observer found in the nonlinear output feedback control
literature (see, e.g.,
\cite{Tor92,EsfKha92,KhaEsf93,TeePra94,LinSab95,MahKha96,Jan96,MahKha97,Isi97,
AtaKha99}). Our observer, however, has the advantage of avoiding the
knowledge of the inverse of the mapping $\HH$, as well as working in
$x$ coordinates, directly. 
\end{rem}
\section{Output Feedback Stabilizing Control}
Consider system (\ref{eq:extended}), by using assumption A2 and Remark
\ref{rem:stabilizability}
we conclude that there exists a smooth stabilizing
control $v=\phi(x,z)=\phi(\chi)$ which makes the origin of
(\ref{eq:extended}) an asymptotically stable equilibrium point with
domain of attraction $\DD$. By the converse Lyapunov
theorem found in \cite{Kur56}, there exists a continuously
differentiable function $V$ defined on $\DD$ satisfying, for all $\chi
\in \DD$, 
\begin{gather}
\alpha_1 ( \| \chi \| ) \leq V(\chi) \leq \alpha_2(\| \chi \|)
\label{eq:V:1}\\
\lim_{\D \chi \tends \partial \DD} \alpha_1 (\| \chi \|) = \infty \label{eq:V:1'}\\
\DER{V}{\chi} \left( f_e (\chi) + g_e \, v\right)  \leq -\alpha_3(\|\chi\|) \label{eq:V:2}
\end{gather}
where $\alpha_i , \, \,i = 1, 2, 3$ are class $\K$
functions (see \cite{Kha96_1} for a definition), and $\partial \DD$
stands for the boundary of the set $\DD$.
Define compact sets  
$\Omega_{c_1}$, $\Omega_{c_2}$, $\Omega_{c_2}^x$, and $\Omega_{c_2}^z$  as follows
%
\[
\Omega_{c_1} \eqdef  \{ \chi  \,|\, V
\leq c_1 \}  , \, \Omega_{c_2} \eqdef \{\chi \,|\, V  \leq c_2 \},  \,
\Omega_{c_2}^x  \eqdef\{ x \in \Re^n \,|\, \chi \in \Omega_{c_2}
\}, \, \Omega_{c_2}^z  \eqdef\{ z \in \Re^{n_u}  \,|\, \chi \in
\Omega_{c_2} \}
\]
where $ c_2 > c_1 > 0$. Next, the following assumption is needed.
\ass{3}. Assume $c_2$ can be selected   so that the following
conditions are satisfied:
\begin{equation*}
\begin{array}{ll}
1. & \HH(\Omega_{c_2}^x , z) \subset C_\xi(z) \subset  
\HH(\X , z), \text{ for all }  z \in \Omega_{c_2}^z, \text{ for some
convex compact  $C_\xi(z)$}
\label{eq:projection_set}\\
2. & \Omega_{c_2}^z \subset \U 
\end{array}
\end{equation*}
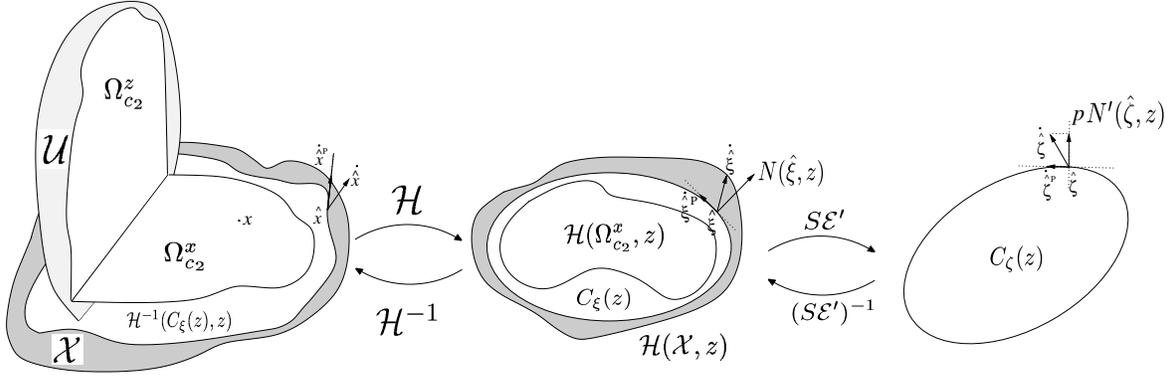
\scalefig{projection}{.95}{The mechanism behind the observer estimates projection. }
\vspace*{-0.5cm}
\begin{rem} See Figure \ref{projection} for a pictorial
representation of the sets in A3. This assumption  represents  
a  basic requirement for output
feedback control. It is satisfied when there exists a sphere
of dimension $n+n_u$,
centered at the origin, which is contained in  $\X \times \U$; this
requirement is satisfied in most practical examples. On the other hand,
Assumption A3 fails when, for example,  the origin belongs to the
boundary  of $\X \times \U$, and thus there is no neighborhood
centered at the origin and contained in $\X \times \U$.
\end{rem}
\subsection{Observer Estimates Projection}
As we already pointed out in Remark 4, in order to isolate the peaking
phenomenon from the system states, the approach generally adopted in
several papers is to  saturate the
control input to prevent it from growing above a given threshold. This
technique, however, does not eliminate the peak in the observer
estimate and,
hence, cannot be used to control general systems like the ones
satisfying assumption A1, since even when the system state lies in the
observable region $\X \times \U \subset \Re^n \times \Re^{n_u}$,
the observer estimates may enter the
unobservable domain where (\ref{eq:observer}) is not well defined. It
appears  that in order to deal with systems that are not completely
observable,
 one has to eliminate the peaking from the observer by
guaranteeing its estimates to be confined in a prespecified compact
set contained in $\X$. 

A very common procedure used in the adaptive control literature (see
\cite{IoaSun95}) to confine vectors of parameter estimates within a
desired convex set is gradient projection. This idea cannot be
directly applied to our problem, mainly because $\dot {\hat x}$ is not 
proportional to the gradient of the observer Lyapunov function and,
thus, the projection cannot be guaranteed to preserve the convergence
properties of the estimate. Inspired by this idea, however, we propose a way to
modify the  $\dot {\hat x}$ equation which confines $\hat x$ to  within a prespecified
compact set while preserving its convergence properties.

Recall the coordinate transformation defined in (\ref{eq:tsf:xi}) and let
\begin{equation}
\xi = \HH(x,z), \s \hat \xi =  \HH(\hat
x,z), \s \tilde \xi = \hat \xi - \xi 
\label{eq:tsf}
\end{equation}
Next, {\em
project}\,\footnote{The projection defined in (\ref{eq:projection}) is
discontinuous in the variable $\hat \xi$, therefore raising the issue
of the existence and uniqueness of its solutions.  We refer the reader
to Remark \ref{rem:projection:smooth}, were this issue is  addressed
and a solution is proposed.}
the observer estimate as follows

\begin{equation}
\begin{split}
&\xPdot = \left[ \DER{\HH}{\hat x} \right]^{-1} \left\{ 
\PP \left(\hat \xi, \dot{\hat \xi}, z, \dot z \right) - \DER{\HH}{z}
\dot z \right\} \\
&\PP(\hat \xi,\dot{\hat \xi},z, \dot z) = 
\begin{cases}
\D \dot{\hat \xi} - \Gamma \frac{N(\hat
\xi) \left(
N(\hat \xi,z)^\top \dot{\hat \xi} + N_z(\hat \xi, z)^\top \dot z
\right)} {N(\hat \xi,z)^\top \Gamma N(\hat \xi,z)} 
& \text{ if }  N(\hat \xi,z)^\top \dot{\hat \xi} + N_z(\hat \xi,
z)^\top \dot z \geq 0 \text{ and }
\hat \xi \in \partial C_\xi(z)\\
\dot{\hat \xi} & \text{ otherwise}
\end{cases}
\end{split}
\label{eq:projection}
\end{equation}
where $\Gamma = ( S \E' )^{-1} ( S \E' )^{-1} $, 
$S = S^\top$ denotes the matrix square root of $P$ (defined in
(\ref{eq:lyap})) and $N(\hat \xi,z)$,  $N_z(\hat \xi,z)$ are the normal vectors to the
boundary of $C_\xi (z)$ \wrt $\xi$ and  $z$, respectively.
The following lemma shows that (\ref{eq:projection}) guarantees
boundedness and preserves convergence for $\hat x$.
\begin{lemma}: If A3 holds and (\ref{eq:projection}) is used:
\begin{itemize}
\item[(i)] Boundedness: $\hat x^P (t) \in \HH^{-1} (C_\xi(z) , z) \subset
\X$ for all $t$, and for all $z \in \Omega_{c_2}^z$.
\end{itemize}
If, in addition, $x \in \Omega_{c_2}^x$ and the assumptions of Theorem
\ref{thm:observer}  are satisfied, then the following is also true
\begin{itemize}
\item[(ii)] Preservation of original convergence characteristics:
properties (i) and (ii) established by Theorem  \ref{thm:observer} remain valid for $\hat x^P$.
\end{itemize}
\label{lem:projection}
\end{lemma}
\pf
In order to prove part (i) we need another coordinate transformation,
$\zeta = S \E' \xi$, (similarly, let $\hat \zeta = S \E' \hat \xi$,
$\tilde \zeta = S \E' \tilde \xi$). Denote by $C_\zeta(z)$ the compact
convex set
in $\zeta$ coordinates, \ie $C_\zeta(z) \eqdef
\left\{\zeta \in \Re^n \,|\, (\E')^{-1} S^{-1} \zeta \in C_\xi(z)
\right\}$, and let $N'(\hat \zeta,z)$, $N'_z(\hat \zeta,z)$ be the normal vectors to the
boundary of $C_\zeta(z)$ \wrt $\zeta$ and $z$, respectively
($N'(\hat \zeta,z)$ is shown in Figure \ref{projection}).
Hence, it is sufficient to show that the
projection (\ref{eq:projection}) will keep $\hat \zeta$ inside
$C_\zeta(z)$, which in turn guarantees that  $\hat x = \HH^{-1} \left
( \hat \xi, z \right)$ is contained in the compact set
$\HH^{-1} \left( C_\xi(z),z \right) \subset \X$. After  coordinate
transformation  (\ref{eq:tsf}) we have that
\begin{align}
\xiPdot &= \frac{d}{d t} \left\{\HH(\hat x^P,z)
\right\} = \left[ \DER{\HH}{\hat x} \xPdot +  \DER{\HH}{ z}
\dot z \right] \\
&= \PP ( \hat \xi, \dot{\hat \xi} )
\end{align}
 In order to relate $N'(\hat
\zeta,z)$, $N'_z(\hat \zeta,z)$  to $N(\hat \xi,z)$, $N'_z(\hat
\zeta,z)$,  notice that the boundary of $C_\xi(z)$ can be
expressed as the set $\partial C_\xi(z)=\{ \xi \in \Re^n \,|\, b(\xi,z) = 0\}$ 
for some continuous function  $b(\xi,z)$, and hence $N(\hat \xi,z) =
\nabla_\xi b |_{(\xi = \hat \xi,z)}$, $N_z(\hat \xi,z) =
\nabla_z b |_{(\xi = \hat \xi,z)}$. Similarly, the boundary of $C_\zeta(z)$ is
the set $\partial C_\zeta(z)=\{ \zeta \in \Re^n \,|\, b( (S \E')^{-1}
\zeta, z) =0 \}$ and
$N'(\hat \zeta,z) = (S \E')^{-1} \nabla_\zeta b|_{(\xi = \hat \xi,z)} = (S
\E')^{-1} N(\hat \xi,z)$, $N'_z(\hat \zeta,z) = N_z(\hat \xi,z)$.  
The expression of the projection
(\ref{eq:projection}) in $\zeta$ coordinates is found by noting that
\begin{equation}
\zetaPdot = S \E' \xiPdot = \left\{ \begin{array}{ll} \D
S \E' \dot{\hat \xi} - (S \E')^{-1} \frac{N \left(
 N^\top \dot{\hat \xi} + N_z^\top \dot z \right) }{N^\top \Gamma N} &
\text{ if }  N^\top \dot{\hat \xi} + N_z^\top \dot z \geq 0 \text{ and }
\hat \xi \in \partial C_\xi(z) \\
S \E' \dot{\hat \xi} & \text{ otherwise}
\end{array}\right. 
\end{equation}
and then substituting $N' = (S \E')^{-1} N$, $N'_z = N_z$,  and
$\dot{\hat \xi} = (S \E')^{-1} \dot{\hat \zeta}$,  to find
that
\begin{equation}
\zetaPdot = \left\{ \begin{array}{ll} \D
\dot{\hat \zeta} - \frac{N' \left(
 N'{}^\top \dot{\hat \zeta} + {N'_z}{}^\top \dot z \right)}{{N'}{}^\top  N'}  &
\text{ if } N'{}^\top \dot{\hat \zeta} + N'_z{}^\top \dot z
 \geq 0 \text{ and }
\hat \zeta \in \partial C_\zeta(z) \\
\dot{\hat \zeta} & \text{ otherwise}
\end{array}\right. 
\label{eq:projection:zeta}
\end{equation}
Next, we show that the domain $C_\zeta(z) = \{ \zeta \in \Re^n \,|\,
b((S\E')^{-1} \zeta,z) \leq 0 \}$ is invariant for
\refeq{eq:projection:zeta}. In order to do that, consider the
candidate Lyapunov function $V_{C_\zeta} = \max\{b\left((S\E')^{-1} \zeta,z\right) , 0\}$
which is positive definite \wrt the set $C_\zeta(z)$, and calculate
its time derivative along the trajectory of
\refeq{eq:projection:zeta} when $\zetaP \in \partial C_\zeta(z)$,
\begin{align}
\dot V_{C_\zeta} &= N'(\hat \zeta,z)^\top \zetaPdot +
N'_z(\zetaP,z)\dot z \\
&= N'{}^\top \dot {\hat \zeta} - \frac{N'{}^\top N' \left( N'{}^\top
\dot{\hat \zeta} + N'_z{}^\top \dot z \right)}{N'{}^\top N'} + N'_z
\dot z \\ 
&= 0  
\end{align}
thus showing that $C_\zeta(z)$ is an invariant set for
\refeq{eq:projection:zeta} and, hence, that
 $\zetaP (t) \in
C_\zeta(z)$ for all $t$, which in turn implies that $\xiP (t) \in
C_\xi(z)$  for all $t$ and, finally, $\xP (t) \in \HH^{-1} (C_\xi(z),z)$ for
all $t$, thus proving part (i) of the theorem.


The proof of part (ii) is based on the
knowledge of a Lyapunov function for the observer in $\tilde \nu$
coordinates (see (\ref{eq:tsf:nu})). Notice   that $\tilde
\zeta = S \tilde \nu$, and $V_o = \tilde \nu^\top P \tilde \nu =
(\tilde \nu^\top S) (S \tilde \nu) = \tilde \zeta^\top \tilde \zeta$. 
We want to show that, in $\tilde \zeta$ coordinates, $\dot V_o<0$
and $\zetaPdot = S \E' \xiPdot$ implies that
$\dot V_o^P \leq \dot V_o$, where $\tilde \zeta^P = \zetaP - \zeta$,
and $V_o^P = \left.{\tilde \zeta^P}\right.^\top  \tilde \zeta^P$.

From (\ref{eq:projection:zeta}), when $\hat \zeta$ is in the
interior of $C_\zeta(z)$, or $\hat \zeta$ is on the boundary of $C_\zeta(z)$
and $N'{}^\top \dot{\hat \zeta} + {N'_z}{}^\top \dot z <0$ (\ie the update is
pointed to the interior of $C_\zeta(z)$),  we have that $V_o =
V_o^P$. Let us consider all the remaining  cases and, since $\hat
\zeta^P = \hat \zeta$ and $\tilde \zeta^P = \tilde \zeta$ (projection
only operates on $\dot {\hat \zeta}$), we have
\begin{equation}
\dot V_o^P = 2 \tilde \zeta^P \tzetaPdot = 2 \tilde \zeta^\top \tzetaPdot
 =  2
\tilde \zeta^\top \left[ \dot {\hat \zeta} - \dot \zeta - p( \hat \zeta,
\dot{\hat \zeta},z,\dot z ) N' (\hat \zeta) \right]
\end{equation}
where $\D p( \hat \zeta,\dot{\hat \zeta}, z, \dot z ) =
\frac{N'{}^\top \dot{\hat \zeta}  + {N'_z}{}^\top \dot z }{N'(\hat
\zeta)^\top  N'(\hat \zeta)}$ is nonnegative since, by
assumption, $N'{}^\top \dot{\hat \zeta}  + {N'_z}{}^\top \dot z  \geq 0$.
Thus,
\begin{equation}
\dot V_o^P = \dot V_o - 2 p \tilde \zeta^\top N'
\end{equation}
Since, by assumption, $x \in \Omega_{c_2}^x$, we have that $\zeta \in
C_\zeta(z)$ and, since  $\hat \zeta$ lies on the boundary of $C_\zeta(z)$,
the vector $\hat \zeta - \zeta$ points outside of
$C_\zeta(z)$ which, by the convexity of $C_\zeta(z)$, implies that $\tilde
\zeta^\top N' \geq 0$, thus concluding the proof of part (ii).
\BBOX
\begin{rem} 
\label{rem:projection:smooth}
In order to avoid the discontinuity in the right hand side of
(\ref{eq:projection}) introduced by $\PP$, one can define $\PP$ to be
the  smooth projection introduced in \cite{PomPra92}. In this case, part
(i) of Lemma \ref{lem:projection} would have to be modified to
\begin{equation}
\hat x^P (t) \in \HH^{-1} (\bar C_\xi(z)
, z), \, \forall z \in \Re^{n_u}
\end{equation}
  where $\bar C_\xi(z) \supset C_\xi(z)$
is a convex set which can be made arbitrarily close to $C_\xi(z)$, and
condition 1 of A3 would have to be replaced by the 
following
\begin{equation}
\begin{array}{ll}
1'.&\HH(\Omega_{c_2}^x,z) \subset C_\xi(z)
\subset  \bar C_\xi(z) \subset \HH(\X , z), \text{ for all }  z \in
\Omega_{c_2}^z,  \text{ for some
convex  compact  sets $C_\xi(z)$ and $\bar C_\xi(z)$}
\end{array}
\end{equation}
\end{rem}
\begin{rem} From the proof of Lemma \ref{lem:projection} we conclude that
(\ref{eq:projection}) performs a   projection for $\hat x$ over the
compact set $\HH^{-1} (C_\xi(z),z)$ which, in general,
is unknown since we do not know $\HH^{-1}$, and is generally {\em not} convex 
(see Figure \ref{projection}). 
It is interesting to note that applying a standard gradient projection for
$\hat x$ over  an arbitrary convex domain does not necessarily
preserve the convergence result (ii) in Theorem \ref{thm:observer}. 
\end{rem}
\subsection{Closed-Loop Stability}
To perform output feedback control we replace the state feedback law
$v=\phi(x,z)$ with $\hat v = \phi (\xP,z)$ which, by the
smoothness of $\phi$ and the fact that $\xP$ is guaranteed to
belong to the compact set $\HH^{-1} (C_\xi(z),z)$, is bounded
provided that $z$ is confined to within a compact set.
Furthermore, the limit (\ref{eq:V:1'}) guarantees that for any compact set $\DD'$
contained in the region of attraction $\DD$, one can choose 
a large enough  $c_1$  so that 
$
\DD' \subset \Omega_{c_1} \subset \Omega_{c_2} \subset \DD
$.
When the observability assumption A1 is
satisfied globally, one can choose any compact $\DD' \subset \DD$; hence, if
A1 and A2 hold globally, $\DD'$ can be any compact set in $\Re^n \times
\Re^{n_u}$.
Taking in account the restriction
 on $c_2$ imposed by Assumption A3,
choose  $\DD'$  to be an arbitrary compact set contained in 
$ \Omega_{c_1}$. 
In the following we will show that $\hat v$ makes the origin of
(\ref{eq:extended}) asymptotically stable and that $\DD'$ is contained
in its region of attraction. The proof is divided in three steps:
\begin{enumerate}
\item (Lemma \ref{lem:invariance}). 
{\em Invariance of $\Omega_{c_2}$ and ultimate boundedness:}
Using the arbitrarily fast rate of convergence of the observer (see
part (ii) in Theorem \ref{thm:observer}), 
we show that any trajectory originating in $\Omega_{c_1}$ cannot exit
the set $\Omega_{c_2}$ and converges in finite time to an arbitrarily
small neighborhood of the origin. Here, Lemma \ref{lem:projection}
plays an important role, in that it guarantees $\hat v$ to be bounded,
and thus it allows us to use the same idea found in \cite{EsfKha92} to
prove stability.
\item (Lemma \ref{lem:asymptotic}). {\em Asymptotic stability of the origin:} 
By using Lemma \ref{lem:invariance} and the exponential stability of
the observer estimate, we prove that the origin of the closed-loop
system is asymptotically stable.
\item (Theorem \ref{thm:closed_loop}). {\em Closed-loop stability:}
Finally, by putting together the results of Lemma \ref{lem:invariance}
and \ref{lem:asymptotic}, we conclude the closed-loop stability proof.
\end{enumerate}
Let $\tilde x \eqdef \hat x - x$, and $\tilde x^P \eqdef \hat x^P -x$, 
and note that part (ii) of Lemma \ref{lem:projection} applies and 
$\tilde  x^P \tends 0$ as $t
\tends \infty$ with arbitrarily fast rate, as long as 
$\chi \in \Omega_{c_2}$.
Moreover, the smoothness of the control law implies that 
\begin{equation}
\| \phi(\hat x,z) - \phi(x,z) \| \leq \bar \gamma \| \tilde x \|
\label{eq:lipschtiz_control}
\end{equation}
for all $x, \hat x \in \HH^{-1}(C_\xi(z),z)$, for all $z \in
\Omega_{c_2}^z$
and some $\bar \gamma > 0$. 
Assume that $\hat x(0) \in \Omega_{c_2}^x$, and let 
$A $ be a positive
constant satisfying $\| \partial V / \partial \chi \| \leq A $ for
all $\chi$ in $\Omega_{c_2}$ (its existence is guaranteed by $V$ being
continuously differentiable).
Now, we can state the following lemma.
\begin{lemma} Suppose that the initial condition $\chi(0)$ is contained in 
$\DD' \subset \Omega_{c_1}$, define the set $\Omega_\epsilon \eqdef \{ \chi :
V(\chi) \leq  d_\epsilon\}$, where $d_\epsilon = 
\alpha_2 \circ \alpha_3^{-1} (\mu \, A \, \bar \gamma
\, \epsilon)$, and choose $\epsilon > 0$ and $\mu >1$ such that
$d_\epsilon < c_1$. Then, there exists a positive scalar $\rho^*, 0<
\rho^* \leq 1$, such
that, for all $\rho \in (0, \rho^*]$, the closed-loop system
trajectories remain confined in $\Omega_{c_2}$,  the set
$\Omega_\epsilon \subset \Omega_{c_1}$ is positively invariant, and is
reached in finite time.
\label{lem:invariance}
\end{lemma}
\pf 
Since $V(\chi(0)) \leq c_1 < c_2$, there exists a time
$T_1 >0$ such that $V(\chi(t)) \leq c_2$, for all $ t \in [ 0 ,
T_1)$.
Choose $T_0$ such that $0<T_0 < T_1$. Then, after noticing that, for all $t \in [0, T_1)$,
$z \in \Omega_{c_2}^z$ and $\hat x^P \in \HH^{-1} (C_\xi(z),z)$ and hence
$v$ is bounded, we can apply Theorem \ref{thm:observer}, part (ii), and conclude that
for any positive $\epsilon$ there exists a positive 
$\rho^*, 0<
\rho^* \leq 1$ such that, for all  $\rho \leq \rho^*$, 
$\| \tilde x^P \| \leq \epsilon, \forall t \in [T_0, T_1)$. Hence,
for all $t \in [ T_0,  T_1 )$, we have that $V(\chi(t))
\leq c_2$ and $\| \tilde x^P(t) \| \leq \epsilon$.  
In order for part (ii) in Theorem \ref{thm:observer} to hold for all $t \geq
T_0$, $\chi$ must belong  to $\Omega_{c_2}$ for all $t \geq 0$. So far
we can only guarantee that $\chi \in \Omega_{c_2}$ for all $t \in [0,
T_1)$ and hence the result of Theorem \ref{thm:observer} applies in this time interval,
only.
Next, we will show that $T_1 = \infty$, \ie $\Omega_{c_2}$ is an
invariant set, so that the result of Theorem \ref{thm:observer} will be guaranteed to hold for all $t \geq 0$.

Consider the Lyapunov function candidate $V$ defined in (\ref{eq:V:1})-(\ref{eq:V:2}).  
Taking its derivative with respect to time and using 
(\ref{eq:lipschtiz_control}),
\begin{align*}
\dot V &= \DER{V}{\chi}  \left[ f_e(x,z) + g_e \phi(x,z) \right] +
\DER{V}{\chi} g_e \left[ \phi(\xP,z) - \phi(x,z) \right]\\
& \leq -\alpha_3(\|\chi\|) + \left\| \DER{V}{\chi} \right\| 
\| \phi(\xP,z) - \phi(x,z) \| \\
& \leq -\alpha_3(\|\chi\|) + A \bar \gamma \| \tilde x^P \| \\
& \leq -\alpha_3(\|\chi\|) + A \bar \gamma \epsilon \\
& \leq -\alpha_3 \circ \alpha_2^{-1} (V) + A \bar \gamma \epsilon 
\end{align*}
for all $t \in [T_0, T_1)$.
When $V \geq d_\epsilon$ we have that 
\[ 
\dot V \leq - (\mu - 1) A \bar \gamma \epsilon
\] 
hence $V$ decays linearly, which in turn implies that
$\chi(t) \in \Omega_{c_2}$ and that $\Omega_\epsilon$ is reached in finite
time.

\BBOX
\begin{rem} The use of the projection for the observer estimate
 plays a crucial role in the proof
of Lemma \ref{lem:invariance}. As $\rho$ is made smaller, the observer peak may
grow larger, thus generating a large control input, which in turn
might drive the system states $\chi$ outside of $\Omega_{c_1}$ in
shorter time. The boundedness of the control input  makes sure that
the exit time $T_1$ is independent of $\epsilon$, since the maximum
size of the $\hat v$ will  not depend on $\epsilon$, thus allowing
us to choose $\epsilon$ independently of $T_1$.
\end{rem}
Lemma \ref{lem:invariance} proves  that the all the trajectories starting in
$\Omega_{c_1}$ will remain confined within $\Omega_{c_2}$ and converge
to an arbitrarily small neighborhood of the origin in finite
time. Now, in order to complete the stability analysis, 
it remains to show that the origin of the output feedback closed-loop system is
asymptotically stable, so that if $\Omega_\epsilon$ is small enough
all the closed-loop system trajectories converge to it.
\begin{lemma} There exists a positive scalar $\epsilon^*$ such that 
for all  $\epsilon \in (0, \epsilon^*]$ 
 all the trajectories starting inside the compact set 
$\Delta_\epsilon \eqdef \{ [\chi^\top, \tilde x^\top]^\top \, | \, V \leq
d_\epsilon \mbox{ and } \| \tilde x \| \leq \epsilon \}$
converge asymptotically to the origin.
\label{lem:asymptotic}
\end{lemma}
\pf Without loss of generality, assume $\epsilon$ is small enough so that 
$\hat x \in \HH^{-1} (C_\xi(z),z)$ and, hence, $\tilde x^P = \tilde x$.
From the proof of Theorem \ref{thm:observer} recall that $\tilde x =  \HH^{-1}
(\hat \xi , z ) -   \HH^{-1} ( \xi , z )$.
Using A1 we have that the mapping $\HH^{-1}$ is
locally Lipschitz. Hence, there exists a neighborhood $N_{\tilde \xi}$
of the origin such that $\| \tilde x \| \leq k_0 \| \tilde \xi \|$,
for all $\tilde \xi \in N_{\tilde \xi}$, and for some
positive constant $k_0$, which, by (\ref{eq:observer_rate}),
 implies that the origin of the $\tilde x$ system is
exponentially stable.
By the converse Lyapunov theorem 
we conclude that there exists a Lyapunov function $V_o'(\tilde x)$ and positive
constants
$\bar c_1 , \bar c_2 , \bar c_3$ such that
\begin{gather*}
\bar c_1 \| \tilde x \|^2 \leq V_o' \leq \bar c_2 \| \tilde x \|^2 \\
\dot V_o' \leq -\bar c_3 \| \tilde x \|^2 
\end{gather*}
Define the positive scalar $\epsilon^*$ such that $\| \tilde x \|
\leq \epsilon^*$ implies $\tilde \xi \in N_{\tilde \xi}$ (the
existence of $\epsilon^*$ is a direct consequence of the fact
that $\HH$ is locally Lipschitz). Next, define
the following composite Lyapunov function candidate
\[ 
V_c (\chi , \tilde x ) = V(\chi) + \lambda \sqrt{V_o'(\tilde x)} , \s \lambda >
\frac{2 \sqrt{\bar c_2} \bar \gamma A}{\bar  c_3}
\]
then, 
\begin{align*}
\dot V_c &\leq -\alpha_3(\|\chi\|) + A \bar \gamma \|\tilde x\| -
\frac{\lambda}{2 \sqrt{V_o'(\tilde x)}} \bar c_3 \| \tilde x \|^2 \\
&\leq -\alpha_3(\|\chi\|) - \left( \frac{\bar c_3 \lambda}{2
\sqrt{\bar c_2}} - A \bar \gamma \right) \|\tilde x\| < 0
\end{align*}
where we have used the fact that $[ \chi^\top , \tilde x^\top ]^\top
\in \Delta_\epsilon$ implies that $\chi \in \Omega_{c_2}$ (provided
$\epsilon$ is small enough), and hence
$\left\|\DER{V}{\chi}\right\| \leq A$.
Since $\dot V_c$ is negative definite, all the $[\chi^\top , \tilde x^\top]^\top$
trajectories starting in $\Delta_\epsilon$ will converge asymptotically to
the origin. 
\BBOX
We are now ready to state the following closed-loop stability theorem.
\begin{thm} For  the closed-loop system (\ref{eq:extended}), (\ref{eq:observer}),
(\ref{eq:projection}),  satisfying
assumptions A1, A2, and A3, 
  the control
law $\hat v = \phi(\xP, z)$, guarantees that there exists a scalar
$\rho^*, 0 < \rho^* \leq 1$, such that, for all $\rho \in (0,
\rho^*]$, the set $\DD' \times \Omega_{c_2}^x$ is contained in the
region of attraction of the origin $(x=0, z=0, \hat x=0)$.
\label{thm:closed_loop}
\end{thm}
\pf
By Lemma \ref{lem:asymptotic}, there exists $\epsilon^* > 0$ such that, for all $\epsilon
\in (0, \epsilon^*]$, $\Delta_\epsilon$ is a region of attraction for the
origin. Use Lemma \ref{lem:invariance} and the fact that $\chi(0) \in \DD' \subset \Omega_{c_1}$ 
to find $\rho^*, 0 < \rho^* \leq 1$, so that for all $\rho \in (0, \rho^*]$ the
state trajectories enter $\Delta_\epsilon$ in finite time. This concludes
the proof of the theorem.
\BBOX

\begin{rem} Theorem \ref{thm:closed_loop} proves regional stability of
the closed-loop
system, since given an observability subspace $\X \times \U$, and provided
condition 1 of A3 is satisfied, the  control law
$\hat v$, together with (\ref{eq:observer}) and 
(\ref{eq:projection}), make the compact set 
$\DD' \times \Omega_{c_2}^x$  a basin of attraction for the origin of the closed-loop
system.
The size of the region of attraction $\DD'\times \Omega_{c_2}^x$ depends on
the size of the set $\X \times \U$ (see condition
(\ref{eq:projection_set})). If  A1 is satisfied globally (as in
\cite{Tor92,TeePra94,AtaKha99}), or $\X \times \U$ is large enough,
then Theorem  \ref{thm:closed_loop} guarantees that the domain of attraction $\DD$ of the
closed-loop system under state feedback is recovered by the output
feedback controller, in that $\DD'$ can be chosen to be any arbitrary
compact set contained in $\DD$. If, besides being completely uniformly
observable, system (\ref{eq:system}) is also stabilizable (and, therefore,
A2 holds globally), then $\DD = \Re^n \times \Re^{n_u}$ and the
result of Theorem \ref{thm:closed_loop}  becomes semi-global, in that 
$\Omega_{c_1}$ and $\Omega_{c_2}$ can be chosen arbitrarily large,
thus achieving the same property of the controller found in \cite{TeePra94}.
\end{rem}
\begin{rem} Analogous to the result in \cite{EsfKha92, AtaKha99},
Theorem \ref{thm:closed_loop}  proves a separation principle for nonlinear systems: given a
stabilizing state feedback controller, the performance of the output
feedback controller recovers the one under state feedback provided
that the parameter $\rho$ is chosen small enough.
Furthermore, by slight modification of Theorem 3 in \cite{AtaKha99},
it is easy to show that the closed-loop system trajectories under
output feedback approach the trajectories under state feedback as
$\rho \tends 0$. In conclusion, the output feedback controller
presented here achieves the same recovery properties of the one in
\cite{AtaKha99} for the more general class of SISO incompletely
observable systems.
We must point out, however, that we assume to have perfect knowledge
of the system dynamics, whereas the results  in \cite{AtaKha99} admit
model uncertainties. We opted not to include model
uncertainties to better illustrate the underlying principles of our
 approach;  it is an easy exercise to show that analogous results to
the ones in
\cite{AtaKha99} hold  when the model uncertainty is included.
\end{rem}
\begin{rem} If system (\ref{eq:system}) is completely uniformly
observable (and hence $\X \times \U = \Re^n \times \Re^{n_u}$), one
can omit the observer estimates projection (\ref{eq:projection})
and use  the standard control input saturation in 
\cite{EsfKha92,Jan96, KhaEsf93,MahKha96,
MahKha97,AtaKha99,Isi97,LinSab95} and prove that, after minor modifications in
Lemma \ref{lem:invariance} and \ref{lem:asymptotic},  Theorem
\ref{thm:closed_loop} would still be valid. In this case, the theory
provided in this paper shows how to achieve output feedback
stabilization with a nonlinear observer, thus avoiding the analytical
knowledge of $\HH^{-1}$.
\end{rem}
\section{Example}
Consider the following state \FLE dynamical system:
\begin{equation}
\begin{split}
&\dot x_1 = x_2 \\
&\dot x_2 = (1+x_1) \exp(x_1^2) + u -1 \\
&y=(x_2-1)^2
\end{split}
\label{eq:example}
\end{equation}
The control input appears in the first derivative of the output:
\[
\dot y = 2(x_2-1)(1+x_1)\exp(x_1^2)+2(x_2-1)(u-1)
\]
Notice, however, that 
the coefficient multiplying $u$ vanishes 
when $x_2=1$, and hence  system (\ref{eq:example}) does not
have a well-defined relative degree everywhere.
The observability mapping $\HH$ is given by
\begin{equation}
y_e = \begin{bmatrix}
y \\ \dot y
\end{bmatrix}
= \begin{bmatrix}
(x_2-1)^2 \\
2 (x_2-1) (1+x_1) \exp(x_1^2) + 2(x_2-1) (u-1)
\end{bmatrix}
\label{eq:example:observability}
\end{equation}
The first equation in (\ref{eq:example:observability})
is invertible for all $x_2 < 1$, and its inverse is given by $x_2 = 1 - \sqrt{y}$.
Substituting $x_2$ into the second equation in
(\ref{eq:example:observability}) and isolating the term in $x_1$, we get
\begin{equation}
(1+x_1) \exp(x_1^2) = \frac{\dot y + 2 \sqrt{y} (u-1)}{-2 \sqrt{y}}
\label{eq:not_invertible}
\end{equation}
Since $(1+x_1) \exp(x_1^2)$ is a strictly increasing function, it
follows that (\ref{eq:not_invertible}) is invertible for all $x_1 \in
\Re$, however, an analytical solution to this equation cannot be found.
In conclusion, Assumption A1 is satisfied
on the domain $\X \times \U = \{ x \in \Re^2 \,|\, x_2<1\} \times \Re
$, but an analytical inverse $x = \HH^{-1} (y_e,u)$ is not known.
The fact that system (\ref{eq:example}) is  incompletely observable,
together with the non-existence of an analytical inverse to 
(\ref{eq:example:observability}), prevents the application of  the
output feedback control approaches in 
\cite{EsfKha92,Tor92,KhaEsf93,TeePra94,LinSab95,MahKha96,Jan96,MahKha97,
Isi97,AtaKha99}.

From (\ref{eq:example:observability}) we have that $n_u=1$, therefore
we add one integrator on the input side,
\[ 
\dot z_1 = v, \s u = z_1 \]
The resulting extended system can be linearized by letting $x_3 =
(1+x_1)\exp(x_1^2)+z_1-1$ and rewriting the system in new coordinates
$x_e \eqdef [x_1, x_2, x_3]^\top$:
\begin{equation}
\begin{split}
&\dot x_1 = x_2 \\
&\dot x_2 = x_3 \\
&\dot x_3 = x_2 \exp(x_1^2) (2 x_1^2 + 2 x_1 + 1) + v
\end{split}
\label{eq:example:extended}
\end{equation}
Choose $v = -x_2 \exp(x_1^2) (2 x_1^2 + 2 x_1 + 1)- K x_e$, where $K =
[1, 3, 3]$, so that the closed-loop system becomes $\dot x_e = (A_c -B_c
K) \, x_e$ with poles placed at $-1$.  Then, the origin $x_e=0$ is
a globally asymptotically equilibrium point of
(\ref{eq:example:extended}), and Assumption A2 is 
satisfied with $\DD = \Re^3$. 

Next, we will seek to find a set $C_\xi(z)$ satisfying Assumption
A3. Let $\bar P$ be the solution of the Lyapunov equation associated
to $A_c - B_c K$, so that a Lyapunov function for system
(\ref{eq:example:extended}) is $V = x_e^\top \bar P x_e$, and any set 
$\Omega_c \eqdef \left\{ x_e \in \Re^3 \,|\, V(x_e) \leq c \right\}$,
with $c>0$, is a region of attraction for the origin.
Observe further that $\chi = [x_1, x_2, z_1]^\top = [x_1, x_2, x_3 -
(1+x_1)\exp(x_1^2)+1]^\top$, and choose $c_2 =  \omega^2
\lambda_{\min}(\bar P)$, where $0 < \omega < 1$. This choice of $c_2$
makes sure that $\Omega_{c_2}$ is contained in the sphere $S_\omega
\eqdef \{ \chi \in \Re^3 \,|\, \|x_e\| \leq \omega \}$ inside which we
have $x_2 \leq \omega$, implying that $x \in \X$. 
When $\chi \in S_\omega$ it is easy to show that
\begin{gather*}
(1-\omega)^2 \leq y \leq (1+\omega)^2 \\
-2 \omega (\omega+1) \leq \dot y \leq 2 \omega (\omega +1)
\end{gather*}
and hence Assumption A3 is verified with $c_2 =  \omega^2
\lambda_{\min}(\bar P)$ and 
\[
C_\xi = \left\{ \xi \in \Re^2 \,|\, (1-\omega)^2 \leq \xi_1 \leq
(1+\omega)^2, \,-2 \omega (\omega+1) \leq \xi_2 \leq 2 \omega (\omega
+1) \right\}
\]
Notice that here a set $C_\xi$ is found which is independent of $z$. In
general, we allow $C_\xi$ to depend on $z$ to make condition 1 in A3
less restrictive.
\twofigcap{simu_2}{$\rho=0.2$.}{simu_05}{$\rho=0.05$.}
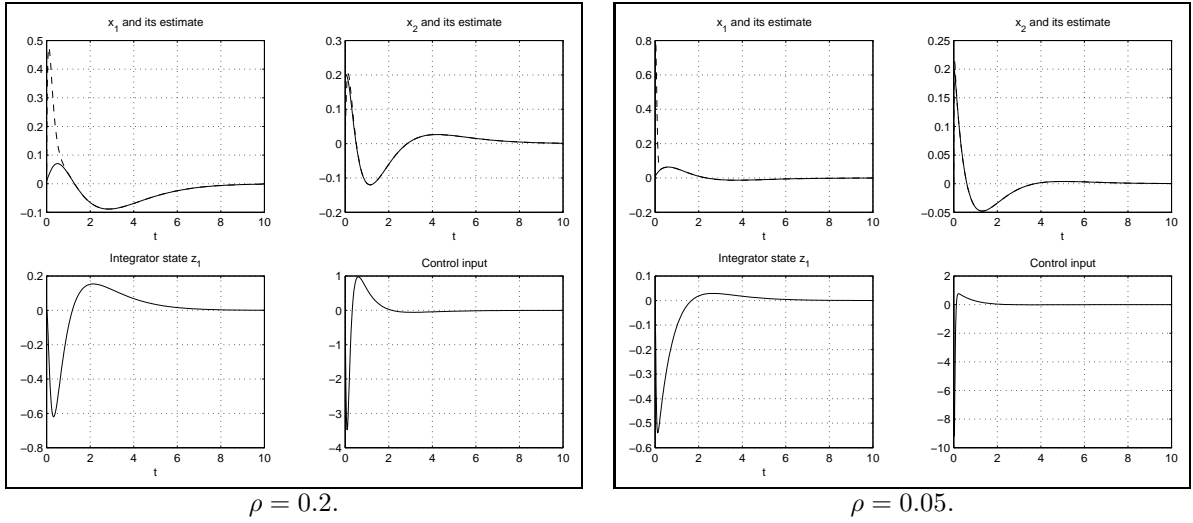
{Closed-loop trajectories under output feedback.}{.45}
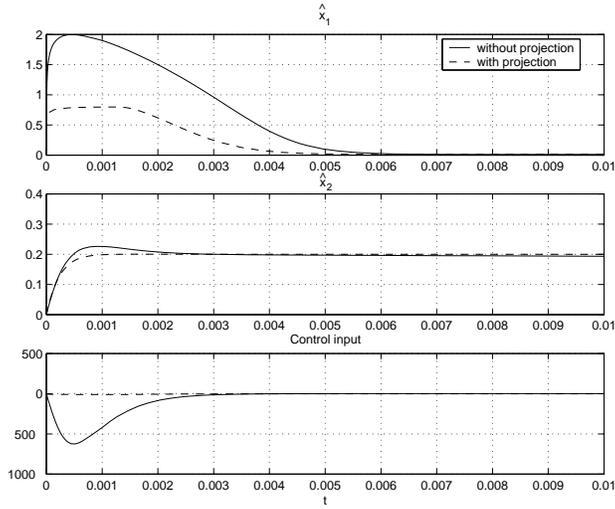
\scalefig{xi}{.5}{Observer states during the
initial peaking phase with and without projection, $\rho = 10^{-3}$.}
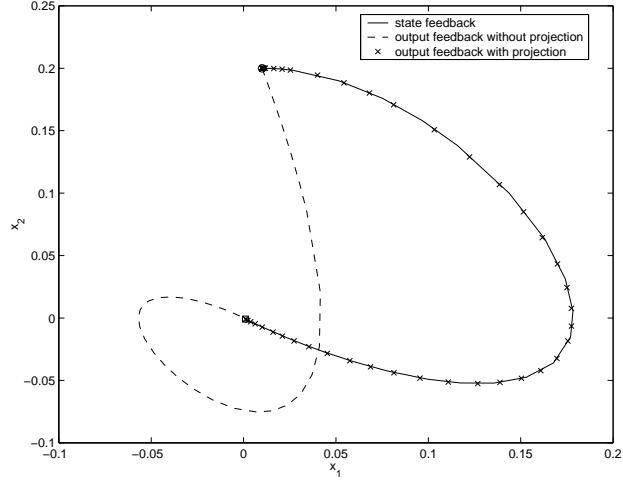
\scalefig{phase_plane}{.5}{System trajectories in the phase plane
with and without projection, $\rho= 10^{-4}$.}
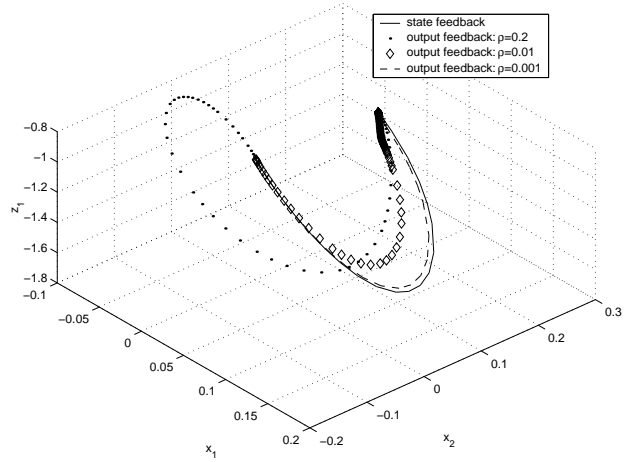
\scalefig{3D}{.5}{Closed-loop trajectories in the three
dimensional space for decreasing values of $\rho$.}

For our simulations we choose $\omega = 0.9$ and the initial condition
of the extended system is set to $x_1(0) = 0.01, x_2(0)= 0.2,
z_1(0)=1.01$, which is contained inside $\Omega_{c_2}$, so that
Theorem \ref{thm:closed_loop} can be applied. Finally,
we choose the observer gain $L$ to be $[4, 4]^\top$, so that its
associated polynomial is Hurwitz with both poles placed at $-2$.
We present four
different situations to illustrate four features of our output
feedback controller:

\begin{enumerate}
\item {\bf Arbitrary fast rate of convergence of the
observer}. 
Figure \ref{simu_2} shows the evolution of the
$\chi$-trajectory, as well as the control input $v$, for $\rho = 0.2$
and $\rho = 0.05$. The convergence in the latter case is faster, as
predicted by Theorem \ref{thm:observer} (see Remark 4).
\item {\bf Observer estimate projection}. 
Figure \ref{xi} shows
the evolution of $\hat x$ and $v$ for $\rho = 10^{-3}$  with and  
without projection. The projection algorithm successfully eliminates  the
peak in the observer states, thus yielding a bounded control input,
as predicted by the result of Lemma \ref{lem:projection}.
\item {\bf Observer estimate projection and closed-loop stability}.
In Figure \ref{phase_plane} a phase plane plot for $x$ is shown
with and without observer projection when $\rho = 10^{-4}$. The small
value of $\rho$ generates a significant peak which, if projection is
not employed, drives the output feedback trajectories away form the 
state feedback ones and, in general, may drive the system to
instability (see Remark 10). On the other hand, using the observer
projection, output feedback and state feedback trajectories are almost indistinguishable.
\item {\bf Trajectory recovery}.
The evolution of the $\chi$-trajectories for decreasing values of
$\rho$, in Figure \ref{3D}, shows that the output feedback
trajectories approach the state feedback ones as $\rho \tends 0$ (see
Remark 12). 
\end{enumerate}

\end{spacing}
\bibliographystyle{ieeetr}
\bibliography{biblio}

\begin{thebibliography}{10}

\bibitem{EsfKha92}
F.~Esfandiari and H.~Khalil, ``Output feedback stabilization of fully
  linearizable systems,'' {\em \IJC}, vol.~56, no.~5, pp.~1007--1037, 1992.

\bibitem{KhaEsf93}
H.~Khalil and F.~Esfandiari, ``Semiglobal stabilization of a class of nonlinear
  systems using output feedback,'' {\em \TAC}, vol.~38, no.~9, pp.~1412--1415,
  1993.

\bibitem{LinSab95}
Z.~Lin and A.~Saberi, ``Robust semi-global stabilization of minimum-phase
  input-output linearizable systems via partial state and output feedback,''
  {\em \TAC}, vol.~40, no.~6, pp.~1029--1041, 1995.

\bibitem{Jan96}
M.~Jankovic, ``Adaptive output feedback control of non-linear feedback
  linearizable systems,'' {\em International Journal of Adaptive Control and
  Signal Processing}, vol.~10, pp.~1--18, 1996.

\bibitem{MahKha96}
N.~Mahmoud and H.~Khalil, ``Asymptotic regulation of minimum phase nonlinear
  systems using output feedback,'' {\em IEEE Trans. on Automatic Control},
  vol.~41, no.~10, pp.~1402--1412, 1996.

\bibitem{MahKha97}
N.~Mahmoud and H.~Khalil, ``Robust control for a nonlinear servomechanism
  problem,'' {\em \IJC}, vol.~66, no.~6, pp.~779--802, 1997.

\bibitem{Isi97}
A.~Isidori, ``A remark on the problem of semiglobal nonlinear output
  regulation,'' {\em \TAC}, vol.~42, no.~12, pp.~1734--1738, 1997.

\bibitem{AtaKha99}
A.~Atassi and H.~Khalil, ``A separation principle for the stabilization of a
  class of nonlinear systems,'' {\em \TAC}, vol.~44, pp.~1672--1687, September
  1999.

\bibitem{Tor92}
A.~Tornamb\`e, ``Output feedback stabilization of a class of non-minimum phase
  nonlinear systems,'' {\em Systems \& Control Letters}, vol.~19, pp.~193--204,
  1992.

\bibitem{TeePra94}
A.~Teel and L.~Praly, ``Global stabilizability and observability imply
  semi-global stabilizability by output feedback,'' {\em Systems \& Control
  Letters}, vol.~22, pp.~313--325, 1994.

\bibitem{Son89}
E.~D. Sontag, ``Remarks on stabilization and input to state stability,'' in
  {\em Proceedings of the IEEE Conference on Decision and Control}, (Tampa,
  FL), pp.~1376--1378, December 1989.

\bibitem{Isi95}
A.~Isidori, {\em Nonlinear Control Systems}.
\newblock London: Springer-Verlag, third~ed., 1995.

\bibitem{KrsKanKok95}
M.~Krsti\'c, I.~Kanellakopoulos, and P.~Kokotovi\'c, {\em Nonlinear and
  Adaptive Control Design}.
\newblock NY: John Wiley \& Sons, Inc., 1995.

\bibitem{CicDalGer93}
G.~Ciccarella, M.~Dalla~Mora, and A.~Germani, ``A {L}uenberger-like observer
  for nonlinear systems,'' {\em \IJC}, vol.~57, no.~3, pp.~537--556, 1993.

\bibitem{Yos66}
T.~Yoshizawa, {\em Stability Theory by Lyapunov's Second Method}.
\newblock The Mathematical Society of Japan, Tokyo, 1966.

\bibitem{Kur56}
J.~Kurzweil, ``On the inversion of {L}japunov's second theorem on stability of
  motion,'' {\em American Mathematical Society Translations, {\rm Series 2}},
  vol.~24, pp.~19--77, 1956.

\bibitem{Kha96_1}
H.~Khalil, {\em Nonlinear Systems}.
\newblock NJ: Prentice-Hall, second~ed., 1996.

\bibitem{IoaSun95}
P.~Ioannou and J.~Sun, {\em Stable and Robust Adaptive Control}.
\newblock Englewood Cliffs, NJ: Prentice-Hall, 1995.

\bibitem{PomPra92}
J.~Pomet and L.~Praly, ``Adaptive nonlinear regulation: Estimation from the
  {L}yapunov equation,'' {\em \TAC}, vol.~37, no.~6, pp.~729--740, 1992.

\end{thebibliography}
\end{document}